\documentclass[amstex,12pt,russian,amssymb]{article}

\usepackage{mathtext}
\usepackage[cp1251]{inputenc}
\usepackage[T2A]{fontenc}
\usepackage[russian]{babel}
\usepackage[dvips]{graphicx}
\usepackage{amsmath}
\usepackage{amssymb}
\usepackage{amsxtra}
\usepackage{latexsym}
\usepackage{ifthen}

\begin{document}

\newcounter{lemma}
\newcommand{\lemma}{\par \refstepcounter{lemma}%
{\bf Лемма \arabic{lemma}.}}

\newcounter{corollary}
\newcommand{\corollary}{\par \refstepcounter{corollary}%
{\bf Следствие \arabic{corollary}.}}

\newcounter{remark}
\newcommand{\remark}{\par \refstepcounter{remark}%
{\bf Замечание \arabic{remark}.}}

\newcounter{theorem}
\newcommand{\theorem}{\par \refstepcounter{theorem}%
{\bf Теорема \arabic{theorem}.}}

\newcounter{proposition}
\newcommand{\proposition}{\par \refstepcounter{proposition}%
{\bf Предложение \arabic{proposition}.}}

\renewcommand{\refname}{\centerline{\bf Список литературы}}

\newcommand{\proof}{{\it Доказательство.\,\,}}

\noindent УДК 517.5

\medskip
{\bf Е.А.~Севостьянов} (Житомирский государственный университет им.\
И.~Франко)

\medskip
{\bf Є.О.~Севостьянов} (Житомирський державний університет ім.\
І.~Франко)

\medskip
{\bf E.A.~Sevost'yanov} (Zhitomir Ivan Franko State University)

{\bf О граничном поведении одного класса отображений на римановых
многообразиях}

{\bf Про граничну поведінку одного класу відображень на ріманових
многовидах}

{\bf On boundary behavior of one class of mappings on Riemannian
manifolds}

Получены теоремы о непрерывном продолжении на границу одного класса
открытых дискретных отображений между произвольными римановыми
многообразиями. В частности, показано, что открытые дискретные
кольцевые $Q$-отображения $f:D\rightarrow D^{\,\prime}$ продолжаются
на границу $\partial D$ области $D,$ как только $\partial D$
локально связна, $\partial D^{\,\prime}$ сильно достижима, а функция
$Q$ имеет конечное среднее колебание в $\partial D.$

Отримано теореми про неперервне продовження на межу одного класу
відкритих дискретних відображень між довільними рімановими
многовидами. Зокрема, показано, що відкриті дискретні кільцеві
$Q$-відображення $f:D\rightarrow D^{\,\prime}$ продовжуються до межі
$\partial D$ області $D,$ як тільки $\partial D$ локально зв'язна,
$\partial D^{\,\prime}$ сильно досяжна, а функція $Q$ має скінченне
середнє коливання в $\partial D.$

Theorems on continuous extension on boundary for one class of open
discrete mappings between Riemannian manifolds are obtained. In
particular, there is proved that, open discrete ring $Q$-mappings
$f:D\rightarrow D^{\,\prime}$ are extend to $\partial D$ whenever
$\partial D$ is locally connected, $\partial D^{\,\prime}$ is
strongly accessible, and a function $Q$ has finite mean oscillation
at $\partial D.$

\newpage
{\bf 1. Введение.} В настоящей статье исследованы некоторые вопросы
граничного поведения на римановых многообразиях отображений, более
общих чем квазиконформные. Здесь изучены открытые дискретные
отображения, в то время как аналогичные проблемы для гомеоморфизмов
исследованы ранее в работе \cite{Af$_1$}. Кроме того, результаты,
аналогичные полученным в настоящей заметке, для случая пространства
${\Bbb R}^n$ опубликованы в \cite{Sev$_1$}.

\medskip
Перейдём к определениям и формулировкам основных результатов.
Следующие понятия могут быть найдены, напр., в \cite{Lee} и
\cite{PSh}. Напомним, что {\it $n$-мерным топологическим
многообразием} ${\Bbb M}^n$ называется хаусдорфово топологическое
пространство со счётной базой, каждая точка которого имеет
окрестность, гомеоморфную некоторому открытому множеству в ${\Bbb
R}^n.$ {\it Картой} на многообразии ${\Bbb M}^n$ будем называть пару
$(U, \varphi),$ где $U$ --- открытое множество пространства ${\Bbb
M}^n$ и $\varphi$ --- соответствующий гомеоморфизм множества $U$ на
открытое множество в ${\Bbb R}^n.$ Если $p\in U$ и
$\varphi(p)=(x^1,\ldots,x^n)\in {\Bbb R}^n,$ то соответствующие
числа $x^1,\ldots,x^n$ называются {\it локальными координатами
точки} $p.$ {\it Гладким многообразием} называется само множество
${\Bbb M}^n$ вместе с соответствующим набором карт $(U_{\alpha},
\varphi_{\alpha}),$ так, что объединение всех $U_{\alpha}$ по
параметру $\alpha$ даёт всё ${\Bbb M}^n$ и, кроме того, отображение,
осуществляющее переход от одной системы локальных координат к
другой, принадлежит классу $C^{\infty}.$

\medskip
Напомним, что {\it римановой метрикой} на гладком многообразии
${\Bbb M}^n$ называется положительно определённое гладкое
симметричное тензорное поле типа $(0,2).$ В частности, компоненты
римановой метрики $g_{kl}$ в различных локальных координатах $(U,
x)$ и $(V, y)$ взаимосвязаны посредством тензорного закона
$$'g_{ij}(x)=g_{kl}(y(x))\frac{\partial y^k}{\partial x^i}
\frac{\partial y^l}{\partial x^j}.$$
{\it Римановым многообразием} будем называть гладкое многообразие
вместе с римановой метрикой на нём. Длину гладкой кривой
$\gamma=\gamma(t),$ $t\in [t_1, t_2],$ соединяющей точки
$\gamma(t_1)=M_1\in {\Bbb M}^n,$ $\gamma(t_2)=M_2\in {\Bbb M}^n$ и
$n$-мерный объём ({\it меру объёма $v$}) множества $A$ на римановом
многообразии определим согласно соотношениям
\begin{equation}\label{eq8}
l(\gamma):=\int\limits_{t_1}^{t_2}\sqrt{g_{ij}(x(t))\frac{dx^i}{dt}\frac{dx^j}{dt}}\,dt,\quad
v(A)=\int\limits_{A}\sqrt{\det g_{ij}}\,dx^1\ldots dx^n\,.
\end{equation}
Ввиду положительной определённости тензора $g=g_{ij}(x)$ имеем:
$\det g_{ij}>0.$ {\it Геодезическим расстоянием} между точками $p_1$
и $p_2\in {\Bbb M}^n$ будем называть наименьшую длину всех
кусочно-гладких кривых в ${\Bbb M}^n,$ соединяющих точки $p_1$ и
$p_2.$ Геодезическое расстояние между точками $p_1$ и $p_2$ будем
обозначать символом $d(p_1, p_2)$ (всюду далее $d$ обозначает
геодезическое расстояние, если не оговорено противное). Так как
риманово многообразие, вообще говоря, не предполагается связным,
расстояние между любыми точками многообразия, вообще говоря, может
быть не определено. Хорошо известно, что любая точка $p$ риманова
многообразия ${\Bbb M}^n$ имеет окрестность $U\ni p$ (называемую
далее {\it нормальной окрестностью точки $p$}) и соответствующее
координатное отображение $\varphi\colon U\rightarrow {\Bbb R}^n,$
так, что геодезические сферы с центром в точке $p$ и радиуса $r,$
лежащие в окрестности $U,$ переходят при отображении $\varphi$ в
евклидовы сферы того же радиуса, а пучок геодезических кривых,
исходящих из точки $p,$ переходит в пучок радиальных отрезков в
${\Bbb R}^n$ (см.~\cite[леммы~5.9 и 6.11]{Lee}, см.\ также
комментарии на стр.~77 здесь же). Локальные координаты
$\varphi(p)=(x^1,\ldots, x^n)$ в этом случае называются {\it
нормальными координатами} точки $p.$ Стоит отметить, что в случае
связного многообразия ${\Bbb M}^n$ открытые множества метрического
пространства $({\Bbb M}^n, d)$ порождают топологию исходного
топологического пространства ${\Bbb M}^n$
(см.~\cite[лемма~6.2]{Lee}). Заметим, что в нормальных координатах
всегда тензорная матрица $g_{ij}(x)$ в точке $p$ --- единичная (а в
силу непрерывности $g$ в точках, близких к $p,$ эта матрица сколь
угодно близка к единичной; см.~\cite[пункт~(c)
предложения~5.11]{Lee}).

\medskip
Всюду далее (если не оговорено противное) ${\Bbb M}^n$ и ${\Bbb
M}_*^n$ -- римановы многообразия с геодезическими расстояниями $d$ и
$d_*,$ соответственно. {\it Кривой} $\gamma$ мы называем непрерывное
отображение отрезка $[a,b]$ (открытого интервала $(a,b),$ либо
полуоткрытого интервала вида $[a,b)$ или $(a,b]$) в ${\Bbb M}^n,$
$\gamma\colon [a,b]\rightarrow {\Bbb M}^n.$ Под семейством кривых
$\Gamma$ подразумевается некоторый фиксированный набор кривых
$\gamma,$ а, если $f\colon{\Bbb M}^n\rightarrow {\Bbb M}_*^n$ ---
произвольное отображение, то
$f(\Gamma)=\left\{f\circ\gamma|\gamma\in\Gamma\right\}.$ Длину
произвольной кривой $\gamma\colon [a, b]\rightarrow {\Bbb M}^n,$
лежащей на многообразии ${\Bbb M}^n,$ можно определить как точную
верхнюю грань сумм $\sum\limits_{i=1}^{n-1} d(\gamma(t_i),
\gamma(t_{i+1}))$ по всевозможным разбиениям $a\leqslant
t_1\leqslant\ldots\leqslant t_n\leqslant b.$ Следующие определения в
случае пространства ${\Bbb R}^n$ могут быть найдены, напр., в
\cite[разд.~1--6, гл.~I]{Va}, см.\ также \cite[гл.~I]{Fu}. Борелева
функция $\rho\colon {\Bbb M}^n\,\rightarrow [0,\infty]$ называется
{\it допустимой} для семейства $\Gamma$ кривых $\gamma$ в ${\Bbb
M}^n,$ если линейный интеграл по натуральному параметру $s$ каждой
(локально спрямляемой) кривой $\gamma\in \Gamma$ от функции $\rho$
удовлетворяет условию $\int\limits_{0}^{l(\gamma)}\rho
(\gamma(s))ds\geqslant 1.$ В этом случае мы пишем:
$\rho\in\mathrm{adm}\,\Gamma.$ {\it Мо\-ду\-лем} семейства кривых
$\Gamma $ называется величина
$$M(\Gamma)=\inf\limits_{\rho\in\mathrm{adm}\,\Gamma}
\int\limits_D \rho ^n (x)\,dv(x).$$
(Здесь и далее $v$ означает меру объёма, определённую в
(\ref{eq8})). При этом, если $\mathrm{adm}\,\Gamma=\varnothing,$ то
полагаем: $M(\Gamma)=\infty$ (см.~\cite[разд.~6 на с.~16]{Va} либо
\cite[с.~176]{Fu}). Свойства модуля в некоторой мере аналогичны
свойствам меры Лебега $m$ в ${\Bbb R}^n.$ Именно, модуль пустого
семейства кривых равен нулю, $M(\varnothing)=0,$ обладает свойством
монотонности относительно семейств кривых, %
$ \Gamma_1\subset\Gamma_2\Rightarrow M(\Gamma_1)\leqslant
M(\Gamma_2), $
а также свойством полуаддитивности:
$ M\left(\bigcup\limits_{i=1}^{\infty}\Gamma_i\right)\leqslant
\sum\limits_{i=1}^{\infty}M(\Gamma_i) $
(см.~\cite[теорема~6.2, гл.~I]{Va} в ${\Bbb R}^n$ либо
\cite[теорема~1]{Fu} в случае более общих пространств с мерами).
Говорят, что семейство кривых $\Gamma_1$ \index{минорирование}{\it
минорируется} семейством $\Gamma_2,$ пишем $\Gamma_1\,>\,\Gamma_2,$
если для каждой кривой $\gamma\,\in\,\Gamma_1$ существует подкривая,
которая принадлежит семейству $\Gamma_2.$
В этом случае,
\begin{equation}\label{eq32*A}
\Gamma_1
> \Gamma_2 \quad \Rightarrow \quad M(\Gamma_1)\leqslant M(\Gamma_2)
\end{equation} (см.~\cite[теорема~6.4, гл.~I]{Va} либо
\cite[свойство~(c)]{Fu} в случае более общих пространств с мерами).

\medskip{}
Следующее определение для случая ${\Bbb R}^n$ может быть найдено,
напр., в работе \cite{SS}. Пусть ${\Bbb M}^n$ и ${\Bbb M}_*^n$ ---
римановы многообразия{\em,} $n\geqslant 2,$ $D$ -- область в ${\Bbb
M}^n,$ $x_0\in D,$ $Q\colon D\rightarrow [0,\infty]$ --- измеримая
относительно меры объёма $v$ функция, и число $r_0>0$ таково, что
замкнутый шар $\overline{B(x_0, r_0)}$ лежит в некоторой нормальной
окрестности $U$ точки $x_0.$ Пусть также $0<r_1<r_2<r_0,$
\begin{equation}\label{eq49***}
A=A(r_1,r_2, x_0)=\{x\in {\Bbb M}^n: r_1<d(x, x_0)<r_2\}\,,
\end{equation}
$S_i=S(x_0,r_i),$ $i=1,2,$ --- геодезические сферы с центром в точке
$x_0$ и радиусов $r_1$ и $r_2,$ соответственно, а
$\Gamma\left(S_1,\,S_2,\,A\right)$ обозначает семейство всех кривых,
соединяющих $S_1$ и $S_2$ внутри области $A.$
%
Пусть $D$ -- область в ${\Bbb M}^n,$ $Q:{\Bbb M}^n\rightarrow [0,
\infty]$ -- измеримая по Лебегу функция, $Q(x)\equiv 0$ при всех
$x\not\in D.$ Отображение $f:D\rightarrow \overline{{\Bbb R}^n}$
будем называть {\it кольцевым $Q$-отоб\-ра\-же\-нием в точке $x_0\in
\partial D,$} если для некоторого $r_0=r(x_0)$ и
произвольных сферического кольца $A=A(r_1,r_2,x_0),$ центрированного
в точке $x_0,$ радиусов: $r_1,$ $r_2,$ $0<r_1<r_2< r_0=r(x_0)$ и
любых континуумов $E_1\subset \overline{B(x_0, r_1)}\cap D,$
$E_2\subset \left(\overline{{\Bbb R}^n}\setminus B(x_0,
r_2)\right)\cap D$ отображение $f$ удовлетворяет соотношению
\begin{equation}\label{eq3*!!}
 M\left(f\left(\Gamma\left(E_1,\,E_2,\,D\right)\right)\right)\ \le
\int\limits_{A} Q(x)\cdot \eta^n(|x-x_0|)\ dm(x) \end{equation}
для каждой измеримой функции $\eta :(r_1,r_2)\rightarrow [0,\infty
],$ такой что
\begin{equation}\label{eq28*}
\int\limits_{r_1}^{r_2}\eta(r)dr \ge\ 1\,.
\end{equation}
Отображения типа кольцевых $Q$-отображений были предложены к
изучению О.~Мартио и изучались им совместно с В.~Рязановым,
У.~Сребро и Э.~Якубовым, см.~\cite{MRSY}, см.\ также \cite{BGMV}.

\medskip
Пусть $\left(X,\,d, \mu\right)$ --- произвольное метрическое
пространство, наделённое мерой $\mu$ и $G(x_0, r)=\{x\in X: d(x,
x_0)<r\}.$ Следующее определение может быть найдено, напр., в
\cite[разд.~4]{RSa}. Будем говорить, что интегрируемая в $G(x_0, r)$
функция ${\varphi}\colon D\rightarrow{\Bbb R}$ имеет {\it конечное
среднее колебание} в точке $x_0\in D$, пишем $\varphi\in FMO(x_0),$
если
%
%
%
%
$$\limsup\limits_{\varepsilon\rightarrow 0}\frac{1}{\mu(G(
x_0,\,\varepsilon))}\int\limits_{G(x_0,\,\varepsilon)}
|{\varphi}(x)-\overline{{\varphi}}_{\varepsilon}|\,
d\mu(x)<\infty,$$
%
%
где
$\overline{{\varphi}}_{\varepsilon}=\frac{1} {\mu(G(
x_0,\,\varepsilon))}\int\limits_{G( x_0,\,\varepsilon)}
{\varphi}(x)\, d\mu(x).$

\medskip
Будем говорить, что граница $\partial D$ области $D$ {\it сильно
достижима в точке $x_0\in
\partial D$}, если для любой окрестности $U$ точки $x_0$ найдется
компакт $E\subset D,$ окрестность $V\subset U$ точки $x_0$ и число
$\delta
>0$ такие, что $M(\Gamma(E,F, D))\ge \delta$ для любого континуума  $F$ в $D,$
пересекающего $\partial U$ и $\partial V$  (см., напр., \cite[разд.
3.8]{MRSY}). В качестве одного из наиболее важных результатов
настоящей статьи, приведём следующее утверждение.

\medskip
\begin{theorem}\label{theor4*!} {\sl\,
Пусть ${\Bbb M}^n$ и ${\Bbb M}_*^n$ --- римановы многообразия,
$n\geqslant 2,$ ${\Bbb M}^n$ компактно, $D$ -- область в ${\Bbb
M}^n,$  $f:D\rightarrow {\Bbb M}_n^*$ -- открытое дискретное
кольцевое $Q$-отоб\-ра\-же\-ние в точке $b\in
\partial D,$ $f(D)=D^{\,\prime},$ область $D$ локально связна в
точке $b,$ $C(f,
\partial D)\subset \partial D^{\,\prime}$ и область $D^{\,\prime}$
сильно достижима хотя бы в одной точке $y\in C(f, b).$ Предположим,
что функция $Q$ имеет конечное среднее колебание в точке $b.$ Тогда
$C(f, b)=\{y\}.$}
\end{theorem}

\medskip
{\bf 2. Вспомогательные сведения}. Всюду далее ${\Bbb M}^n$ ---
риманово многообразие при $n\geqslant 2,$ $d$ --- геодезическое
расстояние на ${\Bbb M}^n,$
 \begin{gather*}
B(x_0, r)=\left\{x\in{\Bbb M}^n: d(x, x_0)<r\right\},\quad S(x_0,r)
= \{ x\,\in\,{\Bbb M}^n: d(x, x_0)=r\},
\end{gather*}
${\rm diam}\,A$ --- геодезический диаметр множества $A\subset {\Bbb
M}^n.$ Всюду далее граница $\partial D$ области $D\subset {\Bbb
M}^n$ и замыкание $\overline{D}$ области $D$ понимаются в смысле
геодезического расстояния $d.$ Перед тем, как мы приступим к
изложению вспомогательных результатов и основной части данного
раздела, дадим ещё одно важное определение (см.~\cite[раздел~3,
гл.~II]{Ri}). Пусть $D$ --- область риманового многообразия ${\Bbb
M}^n,$ $n\geqslant 2,$ $f\colon D \rightarrow {\Bbb M}_*^n$ ---
отображение, $\beta\colon [a,\,b)\rightarrow {\Bbb M}_*^n$ ---
некоторая кривая и $x\in\,f^{\,-1}\left(\beta(a)\right).$ Кривая
$\alpha\colon [a,\,c)\rightarrow D,$ $c\leqslant b,$ называется {\it
максимальным поднятием} кривой $\beta$ при отображении $f$ с началом
в точке $x,$ если $(1)\quad \alpha(a)=x;$ $(2)\quad
f\circ\alpha=\beta|_{[a,\,c)};$ $(3)$\quad если
$c<c^{\prime}\leqslant b,$ то не существует кривой
$\alpha^{\prime}\colon [a,\,c^{\prime})\rightarrow D,$ такой что
$\alpha=\alpha^{\prime}|_{[a,\,c)}$ и $f\circ
\alpha=\beta|_{[a,\,c^{\prime})}.$ Имеет место следующее

\medskip
\begin{proposition}\label{pr7}
{\sl Пусть ${\Bbb M}^n$ и ${\Bbb M}_*^n$ --- римановы многообразия,
$n\geqslant 2,$ $D$ --- область в ${\Bbb M}^n,$ $f\colon
D\rightarrow {\Bbb M}^n_*$ --- открытое дискретное отображение,
$\beta\colon [a,\,b)\rightarrow {\Bbb M}_*^n$ --- кривая и точка
$x\in\,f^{-1}\left(\beta(a)\right).$ Тогда кривая $\beta$ имеет
максимальное поднятие при отображении $f$ с началом в точке $x.$ }
\end{proposition}

\medskip
\begin{proof} Зафиксируем точку $x_0\in {\Bbb M}^n$ и рассмотрим $f(x_0)\in {\Bbb
M}_*^n.$ Поскольку точка $f(x_0)$ принадлежит многообразию ${\Bbb
M}_*^n,$ найдётся окрестность $V$ этой точки, гомеоморфная множеству
$\psi(V)\subset {\Bbb R}^n.$ В силу непрерывности отображения $f,$
найдётся окрестность $U$ точки $x_0,$ такая что $f(U)\subset V.$ С
другой стороны, не ограничивая общности, можно считать, что $U$
гомеоморфна открытому множеству $\varphi(U)$ в ${\Bbb R}^n.$ Можно
также считать, что $\varphi(U)$ и $\psi(V)$ являются областями в
${\Bbb R}^n,$ тогда $f^*=\psi\circ f\circ\varphi^{\,-1}$ ---
открытое дискретное отображение между областями $\varphi(U)$ и
$\psi(V)$ в ${\Bbb R}^n.$ Для таких отображений существование
максимальных поднятий локально вытекает из соответствующего
результата Рикмана в $n$-мерном евклидовом пространстве
(см.~\cite[шаг~2 доказательства теоремы~3.2 гл.~II]{Ri}). Отсюда
вытекает локальное существование максимальных поднятий и на
многообразиях. Глобальное существование максимальных поднятий может
быть установлено аналогично доказательству шага 1 указанной выше
теоремы.~$\Box$
\end{proof}

Следующая лемма доказывалась В.И. Рязановым и Р.Р. Салимовым в
работе \cite{RSa} для случая гомеоморфизмов, и представляет собой
основной результат настоящей работы в наиболее общей ситуации.

\medskip
\begin{lemma}\label{lem1} {\sl\, Пусть ${\Bbb M}^n$ и ${\Bbb M}_*^n$ --- римановы
многообразия, $n\geqslant 2,$ ${\Bbb M}^n$ компактно,  $D$ ---
область в ${\Bbb M}^n,$ $f:D\rightarrow {\Bbb M}_*^n$ -- открытое
дискретное кольцевое $Q$-отоб\-ра\-же\-ние в точке $b\in
\partial D,$ $f(D)=D^{\,\prime},$ область $D$ локально связна в
точке $b,$ $C(f,
\partial D)\subset \partial D^{\,\prime}$ и область $D^{\,\prime}$
сильно достижима хотя бы в одной точке $y\in C(f, b).$ Предположим,
что найдётся $\varepsilon_0>0$ и некоторая положительная измеримая
функция $\psi(t),$ $\psi:(0, \varepsilon_0)\rightarrow (0,\infty),$
такая что для всех $\varepsilon\in(0, \varepsilon_0)$
\begin{equation}\label{eq7***}
0<I(\varepsilon,
\varepsilon_0)=\int\limits_{\varepsilon}^{\varepsilon_0}\psi(t)dt <
\infty
\end{equation}
и при $\varepsilon\rightarrow 0$
\begin{equation}\label{eq5***}
\int\limits_{A(\varepsilon, \varepsilon_0, b)}
Q(x)\cdot\psi^{\,n}(|x-b|)
 \ dv(x) =o(I^{n}(\varepsilon, \varepsilon_0))\,,
\end{equation}
где $A:=A(\varepsilon, \varepsilon_0, b)$ определено в
(\ref{eq49***}). Тогда $C(f, b)=\{y\}.$ }
\end{lemma}

\medskip
\begin{proof}
Предположим противное. Тогда найдутся, по крайней мере, две
последовательности $x_i,$ $x_i^{\,\prime}\in D,$ $i=1,2,\ldots,$
такие, что $x_i\rightarrow b,$ $x_i^{\,\prime}\rightarrow b$ при
$i\rightarrow \infty,$ $f(x_i)\rightarrow y,$
$f(x_i^{\,\prime})\rightarrow y^{\,\prime}$ при $i\rightarrow
\infty$ и $y^{\,\prime}\ne y.$ Отметим, что $y$ и $y^{\,\prime}\in
\partial D^{\,\prime},$ поскольку по условию $C(f,
\partial D)\subset \partial D^{\,\prime}.$ По определению сильно достижимой границы в
точке $y\in \partial D^{\,\prime},$ для любой окрестности $U$ этой
точки найдутся компакт $C_0^{\,\prime}\subset D^{\,\prime},$
окрестность $V$ точки $y,$ $V\subset U,$ и число $\delta>0$ такие,
что
\begin{equation}\label{eq1}
M(\Gamma(C_0^{\,\prime}, F, D^{\,\prime}))\ge \delta
>0
\end{equation} для произвольного континуума
$F,$ пересекающего $\partial U$ и $\partial V.$ В силу предположения
$C(f,
\partial D)\subset \partial D^{\,\prime},$ имеем, что для
$C_0:=f^{\,-1}(C_0^{\,\prime})$ выполнено условие $C_0\cap \partial
D=\varnothing.$ Тогда, не ограничивая общности рассуждений, можно
считать, что $C_0\cap\overline{B(b, \varepsilon_0)}=\varnothing.$
Поскольку область $D$ локально связна в точке $b,$ можно соединить
точки $x_i$ и $x_i^{\,\prime}$ кривой $\gamma_i,$ лежащей в
$\overline{B(b, 2^{\,-i})}\cap D.$ Поскольку $f(x_i)\in V$ и
$f(x_i^{\,\prime})\in D\setminus \overline{U}$ при всех достаточно
больших $i\in {\Bbb N},$ найдётся номер $i_0\in {\Bbb N},$ такой,
что согласно (\ref{eq1})
\begin{equation}\label{eq2}
M(\Gamma(C_0^{\,\prime}, f(\gamma_i), D^{\,\prime}))\ge \delta
>0
\end{equation}
при всех $i\ge i_0\in {\Bbb N}.$ Обозначим через $\Gamma_i$
семейство всех полуоткрытых кривых $\beta_i:[a, b)\rightarrow {\Bbb
R}^n$ таких, что $\beta_i(a)\in f(\gamma_i),$ $\beta_i(t)\in
D^{\,\prime}$ при всех $t\in [a, b)$ и, кроме того,
$B_i:=\lim\limits_{t\rightarrow b-0}\beta_i(t)=\in C_0^{\,\prime}.$
Очевидно, что
\begin{equation}\label{eq4}
M(\Gamma_i)=M\left(\Gamma\left(C_0^{\,\prime}, f(\gamma_i),
D^{\,\prime}\right)\right)\,.
\end{equation}
При каждом фиксированном $i\in {\Bbb N},$ $i\ge i_0,$ рассмотрим
семейство $\Gamma_i^{\,\prime}$ максимальных поднятий
$\alpha_i(t):[a, c)\rightarrow D$ семейства $\Gamma_i$ с началом во
множестве $\gamma_i.$ Такое семейство существует и определено
корректно ввиду предложения \ref{pr7}. Заметим, прежде всего, что
никакая кривая $\alpha_i(t)\in \Gamma_i^{\,\prime},$ $\gamma_i:[a,
c)\rightarrow D,$ не может стремиться к границе области $D$ при
$t\rightarrow c-0$ ввиду условия $C(f,
\partial D)\subset \partial D^{\,\prime}.$ Тогда $C(\alpha_i(t), c)\subset
D.$ Заметим также, что $C(\alpha_i(t), c)\ne \varnothing,$ поскольку
пространство ${\Bbb M}^n$ компактно по предположению леммы.
Предположим теперь, что кривая $\alpha_i(t)$ не имеет предела при
$t\rightarrow c-0.$ Покажем, что предельное множество
$C(\alpha_i(t), c)$ есть континуум в $D.$ Действительно,
по условию Кантора в компакте $\overline{\alpha},$ см. \cite[\S\,
41(I), гл. 4, с. 8--9]{Ku},
$$
G=\bigcap\limits_{k\,=\,1}^{\infty}\,\overline{\alpha\left(\left[t_k,\,c\right)\right)}=
\limsup\limits_{k\rightarrow\infty}\alpha\left(\left[t_k,\,c\right)\right)=
\liminf\limits_{k\rightarrow\infty}\alpha\left(\left[t_k,\,c\right)\right)\ne\varnothing
$$
в виду монотонности последовательности  связных множеств
$\alpha([t_k,\,c))$ и, таким образом, $G$ является связным как
пересечение счётного числа убывающих континуумов
$\overline{\alpha([t_k,\,c))}$ по \cite[теорема 5, \S\,47(II)]{Ku}.

Таким образом, $C(\alpha_i(t), c)$ -- континуум в $D.$ Тогда в силу
непрерывности отображения $f,$ получаем, что $f\equiv const$ на
$C(\alpha_i(t), c),$ что противоречит предположению о дискретности
$f.$

Следовательно, $\exists \lim\limits_{t\rightarrow
c-0}\alpha_i(t)=A_i\in D.$ Отметим, что, в этом случае, по
определению максимального поднятия, $c=b.$ Тогда, с одной стороны,
$\lim\limits_{t\rightarrow b-0}\alpha_i(t):=A_i,$ а с другой, в силу
непрерывности отображения $f$ в $D,$
$$f(A_i)=\lim\limits_{t\rightarrow b-0}f(\alpha_i(t))=\lim\limits_{t\rightarrow b-0}
\beta_i(t)=B_i\in C_0^{\,\prime}\,.$$ Отсюда, по определению $C_0,$
следует, что $A_i\in C_0.$ Погрузим компакт $C_0$ в некоторый
континуум $C_1,$ всё ещё полностью лежащий в области $D,$ см. лемму
1 в \cite{Af$_1$}. За счёт уменьшения $\varepsilon_0>0,$ можно снова
считать, что $C_1\cap\overline{B(b, \varepsilon_0)}=\varnothing.$
Заметим, что функция
$$\eta(t)=\left\{
\begin{array}{rr}
\psi(t)/I(2^{-i}, \varepsilon_0), &   t\in (2^{-i},
\varepsilon_0),\\
0,  &  t\in {\Bbb R}\setminus (2^{-i}, \varepsilon_0)\,,
\end{array}
\right. $$ где $I(\varepsilon,
\varepsilon_0):=\int\limits_{\varepsilon}^{\varepsilon_0}\psi(t)dt,$
удовлетворяет условию нормировки вида (\ref{eq28*}) при
$r_1:=2^{-i},$ $r_2:=\varepsilon_0,$ поэтому, в силу определения
кольцевого $Q$-отоб\-ра\-же\-ния в граничной точке, а также ввиду
условий (\ref{eq7***})--(\ref{eq5***}),
\begin{equation}\label{eq11*}
M\left(f\left(\Gamma_i^{\,\prime}\right)\right)\le \Delta(i)\,,
\end{equation}
где $\Delta(i)\rightarrow 0$ при $i\rightarrow \infty.$ Однако,
$\Gamma_i=f(\Gamma_i^{\,\prime}),$ поэтому из (\ref{eq11*}) получим,
что при $i\rightarrow \infty$
\begin{equation}\label{eq3}
M(\Gamma_i)= M\left(f(\Gamma_i^{\,\prime})\right)\le
\Delta(i)\rightarrow 0\,.
\end{equation}
Однако, соотношение (\ref{eq3}) вместе с равенством (\ref{eq4})
противоречат неравенству (\ref{eq2}), что и доказывает лемму.
\end{proof}$\Box$

\medskip
{\bf 3. Основные результаты.} Следующее утверждение вытекает
из~\cite[лемма~4.1]{RSa} и \cite[следствие~5.1]{ARS}.

\medskip
\begin{proposition}\label{pr3}
{\sl В произвольном римановом многообразии ${\Bbb M}^n$ функция
$Q(x)$ класса $FMO$ в точке $x_0$ удовлетворяет
соотношению~\eqref{eq5***}, где
$\psi(t):=\frac{1}{t\log\frac{1}{t}}.$}
\end{proposition}

\medskip
{\it Доказательство теоремы~{\em\ref{theor4*!}} вытекает из
леммы~{\em\ref{lem1}} на основании
предложения~{\em\ref{pr3}}}.~$\Box$

\medskip{}
{\it Элементом площади} гладкой поверхности $H$ на римановом
многообразии ${\Bbb M}^n$ будем называть выражение вида
$d\mathcal{A}=\sqrt{{\rm det}\,g_{\alpha\beta}^*}\,du^1\ldots
du^{n-1},$
где $g_{\alpha\beta}^*$ --- риманова метрика на $H,$ порождённая
исходной римановой метрикой $g_{ij}$ согласно соотношению
\begin{equation}\label{eq5}
g_{\alpha\beta}^*(u)=g_{ij}(x(u))\frac{\partial x^i}{\partial
u^{\alpha}} \frac{\partial x^j}{\partial u^{\beta}}.
\end{equation}
Здесь индексы $\alpha$ и $\beta$ меняются от $1$ до $n-1,$ а $x(u)$
обозначает параметризацию поверхности $H$ такую, что $\nabla_u x\ne
0.$ Справедливо следующее утверждение, обобщающее теорему
\ref{theor4*!}.

\medskip
\begin{theorem}\label{th1}
{\sl\, Пусть ${\Bbb M}^n$ и ${\Bbb M}_*^n$ --- римановы
многообразия, $n\geqslant 2,$ ${\Bbb M}^n$ компактно,  $D$ ---
область в ${\Bbb M}^n,$ $f:D\rightarrow {\Bbb M}_*^n$ -- открытое
дискретное кольцевое $Q$-отоб\-ра\-же\-ние в точке $b\in
\partial D,$ $f(D)=D^{\,\prime},$ область $D$ локально связна в
точке $b,$ $C(f,
\partial D)\subset \partial D^{\,\prime}$ и область $D^{\,\prime}$
сильно достижима хотя бы в одной точке $y\in C(f, b).$ Предположим,
что при некотором $\delta(x_0)>0$ выполняется равенство
\begin{equation}\label{eq3A}
\int\limits_{0}^{\delta(x_0)}\frac{dt}{\left(\int\limits_{S(x_0,
t)}Q(x)\,d\mathcal{A}\right)^{\frac{1}{n-1}}}=\infty.
\end{equation} Тогда $C(f,
b)=\{y\}.$}
\end{theorem}

\begin{proof} Достаточно показать, что условие~\eqref{eq3A} влечёт
выполнение условия~\eqref{eq5***} леммы~\ref{lem1}. Можно считать,
что $B(x_0, \delta(x_0))$ лежит в нормальной окрестности точки
$x_0.$ Рассмотрим функцию
$$
\psi(t)\quad =\quad\left\{
\begin{array}{rr}
\left(\int\limits_{S(x_0,
t)}Q(x)\,d\mathcal{A}\right)^{\frac{1}{1-n}}, & t\in (0, \delta(x_0)),\\
0, & t\not\in (0, \delta(x_0)).
\end{array}
\right.$$
Заметим теперь, что требование вида~\eqref{eq7***} выполняется при
$\varepsilon_0=\delta(x_0)$ и всех достаточно малых $\varepsilon.$
Далее установим неравенство
\begin{equation}\label{eq4A}
\int\limits_{\varepsilon<d(x, x_0)<\delta(x_0)} Q(x)\psi^n(d(x,
x_0))\,dv(x)\leqslant C\cdot\int\limits_{\varepsilon}^{\delta(x_0)}
\left(\int\limits_{S(x_0,
t)}Q(x)\,d\mathcal{A}\right)^{\frac{1}{1-n}}dt
\end{equation}
при некоторой постоянной $C>0.$ Для этого покажем, что к левой части
соотношения \eqref{eq4A} применим аналог теоремы Фубини. Рассмотрим
в окрестности точки $x_0\in S(z_0, r)\subset {\Bbb R}^n$ локальную
систему координат $z^1,\ldots, z^n,$ $n-1$ базисных векторов которой
взаимно ортогональны и лежат в плоскости, касательной к сфере в
точке $x_0,$ а последний базисный вектор перпендикулярен этой
плоскости. Пусть $r, \theta^1,\ldots, \theta^{n-1}$ сферические
координаты точки $x=x(\theta)$ в ${\Bbb R}^n.$ Заметим, что $n-1$
приращений переменных $z^1,\ldots, z^{n-1}$ вдоль сферы при
фиксированном $r$ равны $dz^1=rd\theta^1,\dots,
dz^{n-1}=rd\theta^{n-1},$ а приращение переменной $z^n$ по $r$ равно
$dz^n=dr.$ В таком случае,
$$dv(x)=\sqrt{{\rm det\,}g_{ij}(x)}r^{n-1}\,dr d\theta^1\dots
d\theta^{n-1}.$$
Рассмотрим параметризацию сферы $S(0, r)$ $x=x(\theta),$
$\theta=(\theta^1,\ldots,\theta^{n-1}),$ $\theta_i\in (-\pi, \pi].$
Заметим, что $\frac{\partial x^{\alpha}}{\partial \theta^{\beta}}=1$
при $\alpha=\beta$ и $\frac{\partial x^{\alpha}}{\partial
\theta^{\beta}}=0$ при $\alpha\ne \beta,$ $\alpha,\beta=1,\ldots,
n-1.$ Тогда в обозначениях соотношения~\eqref{eq5} имеем:
$g_{\alpha\beta}^*(\theta)=g_{\alpha\beta}(x(\theta))r^2,$
$$d\mathcal{A}=\sqrt{\det\, g_{\alpha\beta}(x(\theta))}r^{n-1}d{\theta}^1\ldots d{\theta}^{n-1}.$$
Заметим, что
 \begin{equation}\label{eq6A}
\int\limits_{S(x_0,
r)}Q(x)\psi^n(d(x,x_0))\,d\mathcal{A}=\psi^n(r)r^{n-1}\cdot\int\limits_{\Pi}\sqrt{{\rm
det\,}g_{\alpha\beta}(x(\theta))}Q(x(\theta))\,d\theta^{1}\dots
d\theta^{n-1},
 \end{equation}
где $\Pi=(-\pi, \pi]^{n-1}$ --- прямоугольная область изменения
параметров $\theta^1,\ldots,\theta^{n-1}.$ Согласно сказанному выше,
применяя классическую теорему Фубини (см., напр.,~\cite[разд.~8.1,
гл.~III]{Sa}),
 \begin{multline}\label{eq7A}
\int\limits_{\varepsilon<d(x, x_0)<\delta(x_0)} Q(x)\psi^n(d(x,
x_0))\,dv(x)=\\
=\int\limits_{\varepsilon}^{\delta(x_0)}\int\limits_{\Pi} \sqrt{{\rm
det\,}g_{ij}(x)}Q(x)\psi^n(r)r^{n-1}\,d\theta^1\dots
d\theta^{n-1}dr.
 \end{multline}
Поскольку в нормальных координатах тензорная матрица $g_{ij}$  сколь
угодно близка к единичной в окрестности данной точки, то
$C_2\det\,g_{\alpha\beta}(x)\leqslant\det\,g_{ij}(x)\leqslant
C_1\det\,g_{\alpha\beta}(x).$ Учитывая сказанное и сравнивая
\eqref{eq6A} и \eqref{eq7A}, приходим к соотношению~\eqref{eq4A}. Но
тогда также
$$\int\limits_{\varepsilon<d(x, x_0)<\delta(x_0)}
Q(x)\psi^n(d(x, x_0))dv(x)=o(I^n(\varepsilon, \delta(x_0)))$$
ввиду соотношения~\eqref{eq3A}. Утверждение теоремы следует теперь
из леммы~\ref{lem1}.~$\Box$
 \end{proof}

\medskip
\noindent{{\bf Евгений Александрович Севостьянов} \\
Житомирский государственный университет им.\ И.~Франко\\
кафедра математического анализа, ул. Большая Бердичевская, 40 \\
г.~Житомир, Украина, 10 008 \\ тел. +38 066 959 50 34 (моб.),
e-mail: esevostyanov2009@mail.ru}

\end{document}